\documentclass[11pt]{article}
\usepackage{amssymb,amsmath,latexsym}
\usepackage[dvips]{color}
\usepackage{xcolor}
\usepackage{comment}
\usepackage{mathrsfs}
\allowdisplaybreaks
\numberwithin{equation}{section}
\allowdisplaybreaks
\newtheorem{thm}{Theorem}[section]
\newtheorem{asp:mean}{Assumption}[section]
\newtheorem{lemma}[thm]{Lemma}

\newtheorem{corollary}[thm]{Corollary}

\newtheorem{remark}[thm]{Remark}

\begin{document}
\newtheorem{assumption}{Assumption}

\def\ee{\varepsilon}
\def\qed{{\hfill $\Box$ \bigskip}}
\def\MM{{\cal M}}
\def\BB{{\cal B}}
\def\LL{{\cal L}}
\def\FF{{\cal F}}
\def\EE{{\cal E}}
\def\QQ{{\cal Q}}
\def\AA{{\cal A}}

\def\cB{\mbox{${\cal B}$}}
\def\<{\langle}  \def\>{\rangle}

\def\R{{\mathbb R}}
\def\N{{\mathbb N}}
\def\E{{\bf E}}
\def\F{{\bf F}}
\def\H{{\bf H}}
\def\P{{\bf P}}
\def\Q{{\bf Q}}
\def\S{{\bf S}}
\def\J{{\bf J}}
\def\K{{\bf K}}
\def\F{{\bf F}}
\def\A{{\bf A}}
\def\loc{{\bf loc}}
\def\eps{\varepsilon}
\def\semi{{\bf semi}}
\def\wh{\widehat}
\def\pf{\noindent{\bf Proof} }
\def\dim{{\rm dim}}

\title{Coalescence times for  critical Galton-Watson
processes with immigration\footnote{The research of this project is supported by the National Key R\&D Program of China (No. 2020YFA0712902) } }
\author{ Rong-Li Liu$^{a,}$,\footnote{The research of R. Liu is supported in part by NSFC (Grant No. 12271374)}
 \hspace{1mm}\hspace{1mm} Yan-Xia
Ren,$^{b,}$\footnote{The research of Y.-X. Ren is supported in part  by NSFC (Grant Nos. 12071011 and 12231002) 
and The Fundamental Research Funds for the Central Universities, Peking University LMEQF\hspace{1mm} }
\hspace{1mm}\hspace{1mm}  Yingrui Wang$^c$
}
\date
{\footnotesize
$^a$  School of Mathematics and Statistics, Beijing Jiaotong University, Beijing, 100044, P.R.China. E-mail: \qquad\qquad\qquad\qquad\qquad\qquad {\tt rlliu@bjtu.edu.cn}\\
$^b$ LMAM School of Mathematical Sciences \& Center for
Statistical Science, Peking University, Beijing 100871, P. R. China.
E-mail: {\tt yxren@math.pku.edu.cn} \\
$^c$ School of Mathematics and Statistics, Beijing Jiaotong University, Beijing 100044, P. R. China.
 E-mail: {\tt wyr2120873@gmail.com}
 }
 \maketitle

\noindent{\bf Abstract}\\
\noindent Let $X^I_n$ be the  coalescence time of two particles picked at random
from  the $n$th generation of a critical Galton-Watson
process with immigration, and let $A^I_n$ be the coalescence time of the whole population in the $n$th generation.
In this paper, we study the limiting behaviors of $X^I_n$ and $A^I_n$ as $n\to\infty$.

\bigskip
\noindent{\bf Keywords}\hspace{2mm}
critical Galton-Watson  process,  immigration, coalescence times.
\medskip

\noindent{\bf  2010 MR Subject Classification}
primary 60J68; secondary 62E15; 60F10; 60J80

\bigskip
\section{Introduction and Main Results}
Suppose $(Y_n)_{n\geq 0}$ is  a Galton-Watson
process with offspring distribution
$(p_j)_{j\geq 0}$ and
initial size $Y_0=1$.
For  $n\geq 1$,
conditional on $\{Y_n \geq 2\}$,
pick 2 distinct particles uniformly from the $n$-th generation
 and trace their lines
of descent backward in time.
The common nodes in the two lines are called
  the common ancestors of the two particles.
Let $X_n$ denote the generation of their most recent common ancestor, which is called {\em the pairwise coalescence time}.
Next, for $n\geq 1$,
conditional on $\{Y_n \geq 1\}$,
 we trace the  lines  of descent of all particles
 in generation $n$  backward in time.
The common nodes in the $Y_n$ lines of descent are called the
common ancestors of  all the particles in generation $n$.
Define {\em the total coalescence time} $A_n$
as the generation of the most recent common ancestor of all the particles
in generation $n$.  When $m:=\sum^\infty_{n=0}jp_j=1$ (critical case), $p_1<1$ and $\sigma^2:=\sum^\infty_{n=0}j^2p_j-1<\infty$, Athreya \cite{Atha} proved that  for $u\in (0, 1)$,
\begin{equation}\label{num in Ath}
\lim_{n\to\infty}P\left(\frac{X_n}{n}\ge u\big|Y_n\geq 2\right)=E\left[\frac{\sum^{N_u}_{i=1}\eta^2_i}{(\sum^{N_u}_{i=1}\eta_i)^2}\right],
\end{equation}
where $(\eta_i )_{i\geq 1}$ are independent and identically distributed  exponential random
variables with mean $\sigma^2/2$,
and $N_u$ is independent of $(\eta_i)_{i\geq 1}$ and is
a geometric  random variable of parameter $1-u$
(i.e., $P(N_u = k) =(1-u)u^{k-1}, k\geq 1$).
Athreya \cite{Atha} also proved  the following conditional limit result:
$$
\lim_{n\to\infty}P\left(\frac{A_n}{n}> u\big|Y_n\geq 1\right)=1-u,\quad \mbox{ for }u\in (0,1).
$$
The genealogy of branching processes has
been widely studied. Athreya \cite{Athc,Athb}, Durrett \cite{Du},  Zubkov \cite{Z} also investigated the distributional properties of the coalescence times for Galton-Watson
processes.  Kersting \cite{Ker} gave the genealogy structure of branching processes in random environment.  Harris, Johnston and Roberts
\cite{HJR}, Johnston \cite{J} and Le \cite{Le} investigated the coalescent structure of continuous time Galton-Watson
processes. Hong \cite{H} studied the corresponding results for multitype branching processes.

Suppose $(p_j)_{j\geq 0}$ and $(b_j)_{j\geq 0}$ are probability distributions on
the set $\mathbb N$ of nonnegative integers.
Let $(\xi_{n,i}; n\in \mathbb N, i\in\mathbb N)$ be a doubly infinite family of independent random variables with common distribution $(p_j)_{j\geq 0}$, and let $(I_n)_{n\geq 0}$ be a sequence of independent random variables with common distribution $(b_j)_{j\geq 0}$ which are independent of
$(\xi_{n,i}; n\in \mathbb N, i\in\mathbb N)$ as well.  Let $(Z_n)_{n\geq 0}$ be a Galton-Watson process with immigration (GWPI for short) defined by
\begin{eqnarray*}\label{def of GWPI}
Z_0=I_0,\qquad Z_{n+1}=\sum_{i=1}^{Z_n}\xi_{n,i}+I_{n+1},\qquad n=0,1,\ldots.
\end{eqnarray*}
Here $Z_n$ is the  population size
in generation $n$, and $I_n$ is the number of immigrants in generation $n$.
For each $1\leq i\leq Z_n$,  $\xi_{n,i}$ denotes the number of children of the $i$-th particle in generation $n$.
We assume that all the immigrants have different ancestors.
Set $m=E\xi_{0,1}=\sum_{j=0}^\infty jp_j$.
Then $(Z_n)_{n\geq 0}$ is called  supercritical, critical or subcritical
according to $m>1, m=1$ or $m<1$, respectively. GWPI was
first considered
by Heathcote \cite{HCR} in 1965.
Recently, Wang, Li and Yao
\cite{WLY}
found that the pairwise coalescence time $X_n$ for some supercritical GWPI converges in distribution to a $(0,\infty]$-valued random variable as $n\to\infty$.

In this paper, we consider the coalescence times for
critical GWPI $(Z_n)_{n\geq 0}$.
Unlike the case of a Galton-Watson process starting with one particle,
two randomly picked distinct particles (all particles )
from generation $n$ of a GWPI may not have a common ancestor.
Conditional on $\{Z_n>1\}$,
we pick two distinct particles, say $v_1$ and $v_2$, uniformly from
the $n$th generation and trace their lines of descent
backward in time.
Define the
{\em pairwise coalescence time} for GWPI
\begin{equation}
X^I_n=\left\{\begin{array}{ll}|v|,\quad &\mbox{if the most recent common ancestor of } v_1 \mbox{ and }  v_2 \mbox{ is }v,\\
\infty, \quad &\mbox{otherwise},\end{array}\right.
\end{equation}
where $|v|$  is the generation of $v$.
 Note that
 even if $v_1$ and $v_2$
 are descendants of two distinct particles immigrated to the system at the same time,
we  do not say they have a common ancestor.
Similarly, conditional on $\{Z_n>0\}$, define
{\em the total coalescence time} for GWPI
\begin{equation}
A^I_n=\left\{\begin{array}{ll}|v|,\quad &\mbox{if  the most recent common ancestor of all particles  alive at } n \mbox{ is  } v,\\
\infty, \quad &\mbox{otherwise}.\end{array}\right.
\end{equation}
We will study the asymptotic behaviors of
the distribution of $X^I_n$
conditioned on $\{Z_n>1\}$ and
the distribution of $A^I_n$
 conditioned on $\{Z_n>0\}$.
We will explore the effect of the immigrations on the coalescence times.
Throughout this paper we  suppose the following assumption holds.
\begin{assumption}\label{offspring assumption}
$0<p_0+p_1<1$,  $m=1$, $\sigma^2=\sum_j (j^2-1)p_j<\infty$. $b_0<1$ and $\beta=\sum_j jb_j<\infty$.
\end{assumption}

 We use $\langle g,\mu\rangle$ to denote the integral of
a function $g$ with respect to a Radon measure $\mu$ whenever this integral makes sense.
\begin{thm}\label{main limit thm1}
Suppose Assumption \ref{offspring assumption} holds. Let
$\gamma=2\beta/\sigma^2$.  Define
\begin{eqnarray}\label{fraction}
\phi(j, \mu )=E\Big[\dfrac{\sum_{i=1}^{j} \omega^2_i +\langle f^2, \mu\rangle}{(\sum_{i=1}^{j} \omega_i+\langle f, \mu\rangle)^2}\Big],
\end{eqnarray}
where $f(r)=r, r>0$, and
$(\omega_i)_{i\ge 1}$ are independent exponential random variables with parameter
$\frac{2}{\sigma^2}$.
\begin{itemize}
\item[$(1)$] For $0<u<1$,
\[
\lim_{n\to\infty, k/n\to u} P\left(k\leq X^I_n<n\Big|Z_n>1\right)= E\phi(N^I_u, W),
\]
where $N^I_u$ is a negative binomial random variable with law
\begin{eqnarray}\label{nbl}
P(N^I_u=k)
=\dfrac{(-\gamma)(-\gamma-1)\cdots (-\gamma-k+1)}{k!}(1-u)^{\gamma}(-u)^{k},\ \quad k=0,1,2,\ldots,
\end{eqnarray}
with the convention
$\frac{(-\gamma)(-\gamma-1)\cdots (-\gamma-k+1)}{k!}=1$ when $k=0$, $W$ is a Poisson random measure on $(0,\infty)$ with intensity $\frac{\gamma}{r}e^{-\frac{2}{\sigma^2}r}dr$, and $N^I_u$ and $W$ are independent.
\item[$(2)$]
\[
 \lim_{n\to\infty}P\left(X^I_n<\infty\big|Z_n>1\right)=E\left[\dfrac{\langle f^2, W\rangle}{\langle f, W\rangle^2}\right].
\]
\end{itemize}
\end{thm}
Note that $N_u$ in \eqref{num in Ath} for a critical Galton-Watson process
 only takes positive integer values, while
 $N^I_u$ in Theorem \ref{main limit thm1} can take value $0$ with positive probability.
 In the special case $\gamma=1$, the random number
$N^I_u+1$ and $N_u$ have the same distribution.

We conclude from \cite[Theorem $3$]{pakes71b} (see Lemma \ref{gamma limit}) that $Z_n$ diverges to infinity in probability as $n\to\infty$.
Our second result says that
as $n\to\infty$, the probability that all the particles
of generation $n$ have a common ancestor goes to $0$.
\begin{thm}\label{result for taun}
Suppose Assumption \ref{offspring assumption} holds. Then
\[
 \lim_{n\to\infty}P(A_n^I<\infty|Z_n>0)=0.
 \]
\end{thm}
\section{Some preliminary results}
Recall that $(Y_n)_{n\geq 0}$
is a critical Galton-Watson process with offspring distribution $(p_j)_{j\ge 0}$ starting with $Y_0=1$.
The following result was
proved  in \cite{Ath Ney}.
\begin{lemma}\label{basic lemma}
When $m=1$, $p_1<1$, $\sigma^2=\sum_j (j^2-j)p_j<\infty$,
\begin{equation}\label{rate of ext prob with finit variance}
\lim_{n\to\infty} nP(Y_n>0)=\dfrac{2}{\sigma^2},
\end{equation}
and
for any $t>0$,
\begin{equation}\label{limit of distri with finit var}
\lim_{n\to\infty}P\left(\dfrac{Y_n}{n}>t\big|Y_n>0\right)=e^{-\frac{2t}{\sigma^2}}.
\end{equation}
\end{lemma}
The following result for critical GWPI is from \cite[Theorem 3]{pakes71b}.
\begin{lemma}\label{gamma limit}
Suppose Assumption \ref{offspring assumption} holds. Put $\gamma=\frac{2\beta}{\sigma^2}$. Then, as $n\to\infty$,
$\frac{Z_n}{n}$ converges
 in distribution to
 a Gamma random variable with parameters $(2/\sigma^2, \gamma)$, whose density function is
\begin{equation}\label{density of gamma dis}
h(t)=\dfrac{2}{\sigma^2\Gamma(\gamma)}\left(\frac{2t}{\sigma^2}\right)^{\gamma-1}e^{-\frac{2t}{\sigma^2}},\qquad t>0.
\end{equation}
\end{lemma}
The above lemma
implies that $\lim_{n\to\infty}P(Z_n>0)=1.$
The rate that $1-P(Z_n>0)$ converges  to $0$ was investigated in \cite{pakes71a}.

From the construction \eqref{def of GWPI} of the GWPI $(Z_n)_{n\geq 0}$, for any $0\leq k<n$, $Z_n$ can be rewritten as
\begin{equation}\label{decomp. at k}
Z_n=\sum_{i=1}^{Z_k}Y_{n,k,i}+\sum_{j=k+1}^n \sum_{l=1}^{I_j}Y_{n-j, l}^{(j)},
\end{equation}
where $Y_{n,k,i}, i=1,2,\ldots, $ are independent and have the same distribution as $Y_{n-k}$, and
for $0\leq j\leq n$, $Y_{n-j,l}^{(j)},  l=1,2\ldots$,
are independent and have the same distribution as $Y_{n-j}$. 
Note that $Y_{n,k,i}$ represents the number of descendants in generation $n$ of the $i$th particle in generation $k$, and 
$Y_{n-j,l}^{(j)}$ represents the number of descendants in generation $n$ of the $l$th particle in the $I_j$ immigrants in generation $j$. 
For any non-negative integer $m$, set $(m)_2=m(m-1)$. Notice that $(m)_2=0$ when $m=0$ or $1$.  Starting from the representation \eqref{decomp. at k}, the distribution of the pairwise coalescence time
$X_n^I$,
given $\{Z_n>1\}$, has the following expression.
\begin{lemma}\label{main thm 1}
For any $0\leq k<n$,
\[
P\left(k\leq X^I_n<n\big|Z_n>1\right)
=E\Big[\dfrac{\sum_{i=1}^{Z_k}\big(Y_{n,k,i}\big)_2+\sum_{j=1+k}^{n-1}\sum_{l=1}^{I_j}\big(Y_{n-j,l}^{(j)}\big)_2}{(Z_n)_2}\Big|Z_n>1\Big],
\]
with the convention that the second term in the numerator equals $0$ when $k>n-2$.  In particular,
\[
P\big(X_n^I<\infty|Z_n>1\big)
=E\Big[\dfrac{\sum_{j=0}^{n-1}\sum_{l=1}^{I_j}\big(Y_{n-j,l}^{(j)}\big)_2}{(Z_n)_2}\Big|Z_n>1\Big].
\]
\end{lemma}
\noindent {\bf Proof.}
For $0\leq k<n$, the event $\{k\leq X^I_n<n\}$ occurs
if and only if
either the two randomly picked particles from generation $n$ are both descendants of a
particle in the $k$th generation,
 or they are both descendants of a particle immigrated
 into the system between generation $k+1$ and
generation $n-1$.
The number of choices of the two particles from the descendants
of the $i$th particle
in generation $k$
is $(Y_{n,k,i})_2$, and therefore the total number is $\sum_{i=1}^{Z_k}(Y_{n,k,i})_2$ with the convention that the sum is $0$ if $Z_k=0$.  The number of choices of the two particles from the
descendants of the $l$th particle immigrated
 into the system in generation $j$ for $k+1\leq j<n$ is $(Y_{n-j,l}^{(j)})_2$, and the total number is $\sum_{j=k+1}^{n-1}\sum_{l=1}^{I_j}(Y_{n-j,l}^{(j)})_2$.  Also, the total number of choices of the two
particles from the $n$th generation is $(Z_n)_2$.  Thus for any $n\geq 1$ and $0\leq k<n$,
conditional on $\{Z_n>1\}$, the probability of
$\{k\leq X_n^I<n\}$ is given by
\[
P(k\leq X^I_n<n|Z_n>1)
=E\Big[\dfrac{\sum_{i=1}^{Z_k}(Y_{n,k,i})_2+\sum_{j=k+1}^{n-1}\sum_{l=1}^{I_j}(Y_{n-j,l}^{(j)})_2}{(Z_n)_2}\Big|Z_n>1\Big].
\]
Since $Z_0=I_0$, we have
$Y_{n,0,i}=Y_{n,i}^{(0)}$, $i=1,\ldots, I_0$.
 Taking $k=0$ in the above identity, we obtain
\[
P(X^I_n<\infty|Z_n>1)
=P(X^I_n<n|Z_n>1)
=E\Big[\dfrac{\sum_{j=0}^{n-1}\sum_{l=1}^{I_j}(Y_{n-j,l}^{(j)})_2}{(Z_n)_2}\Big|Z_n>1\Big].
\]
\qed

Let $\mathcal M$ be the space of finite measures on $[0,\infty)$ equipped with
the topology of weak convergence.
Let $C_b[0,\infty) (C_b^+[0,\infty))$ be the space of
 bounded continuous (nonnegative bounded continuous) functions on $[0,\infty)$.
Then for any $g\in C_b[0,\infty)$,
the map $\pi_g: \mu\to \langle g,\mu\rangle$ on $\mathcal M$ is continuous.
 For random measures $\eta_n, \eta\in\mathcal M$,  $n=1,2,\ldots$,
 $\eta_n$ converges to $\eta$ in distribution as $n\to\infty$ is equivalent to
$\langle g,\eta_n\rangle\stackrel{d}{\to}\langle g,\eta\rangle$
for all $g\in C_b^+[0,\infty)$.
We refer the readers to \cite[p.109]{Kal} for more details.  Let $\mathcal F_k$ be the $\sigma$-algebra generated by $\xi_{i,j}, i<k, j=1,2,\ldots$, and
$I_j, j=0,1,\ldots, k$.
Then $\mathcal F_k$ contains all information up to generation $k$.  For $k\geq 0$, given ${\mathcal F}_k$,
 $(Y_{n,k,i})_{n\geq k}$, $i=1,2,\ldots$, are
 independent critical Galton-Watson processes with initial value $1$ at
generation $k$.
\begin{lemma}\label{Lemma of crm}
Suppose Assumption \ref{offspring assumption} holds. If  $\frac{k}{n}\to u$ as $n\to\infty$ for some $u\in(0,1)$, then
as $n\to\infty$, the random measure
$$
V_{n,k}(\cdot)=\sum_{i=1}^{Z_k}{\rm I}_{\{Y_{n,k,i}>0\}}\delta_{\frac{Y_{n,k,i}}{n-k}}(\cdot)\in\mathcal M
$$
converges in distribution to the random measure
$V_u:=\sum_{i=1}^{N^I_u}\delta_{\omega_i}(\cdot)\in \mathcal M$
with the convention that $V_u=0$ when
$N^I_u=0$, where $(\omega_i)_{i\ge 1}$
are independent exponential random variables with parameter $\frac{2}{\sigma^2}$, and
$N^I_u\in\mathbb N$ is independent of
$(\omega_i)_{i\geq 1}$ with the law given by \eqref{nbl}.
\end{lemma}
\noindent{\bf Proof.}  Suppose $g\in C_b^+[0,\infty)$.  For any $0\leq k<n$, let
\[
L_{n,k}(g)=\exp\left\{-\langle g,V_{n,k}\rangle\right\}=\exp\Big\{-\sum_{i=1}^{Z_k}g
\big(\frac{Y_{n,k,i}}{n-k}\big){\rm I}_{\{Y_{n,k,i}>0\}}\Big\},
\]
and set $S_{n,k}g=E\Big(\exp\big\{-g\big(\frac{Y_{n-k}}{n-k}\big){\rm I}_{\{Y_{n-k}>0\}}\big\}\Big)$.  Then we have
\begin{eqnarray*}\label{Lap of branching}
E[L_{n,k}(g)|\mathcal F_k]=E[L_{n,k}(g)|Z_k]
=\Big[E\Big(\exp\big\{-g\Big(\frac{Y_{n-k}}{n-k}\Big){\rm I}_{\{Y_{n-k}>0\}}\big\}\Big)\Big]^{Z_k}
=(S_{n,k}g)^{Z_k}.
\end{eqnarray*}
Let $q_n=P(Y_n>0)$ be the survival probability of the process $(Y_k)_{k\geq 0}$ in generation $n$.
Then we have
\begin{eqnarray*}
S_{n,k}g &=& E\Big[\exp\Big\{- g\Big(\frac{Y_{n-k}}{n-k}\Big)\Big\}\Big|Y_{n-k}>0\Big]q_{n-k}+(1-q_{n-k})\\
&=&1-q_{n-k}\Big[1-E\Big(\exp\Big\{- g\Big(\frac{Y_{n-k}}{n-k}\Big)\Big\}\Big|Y_{n-k}>0\Big)\Big].
\end{eqnarray*}
It follows from \eqref{limit of distri with finit var} that for any $g\in C_b^+[0,\infty)$ and $u\in(0,1)$,
\[
\lim_{n\to\infty, k/n\to u}E\Big[\exp\Big\{-g\Big(\frac{Y_{n-k}}{n-k}\Big)\Big\}\Big|Y_{n-k}>0\Big]=\dfrac{2}{\sigma^2}\int_0^\infty e^{-g(r)}e^{-\frac{2r}{\sigma^2}}dr=:L(g).
\]
By the dominated convergence theorem for convergence in distribution, we have that
\begin{eqnarray*}
&&\lim_{n\to\infty, k/n\to u}E[L_{n,k}(g)]=\lim_{n\to\infty, k/n\to u}E\big[E\left(L_{n,k}(g)|\mathcal F_k\right)\big]=\lim_{n\to\infty, k/n\to u} E\big[(S_{n,k}g)^{Z_k}\big]\\
&=&E\lim_{n\to\infty, k/n\to u}\Big[\Big(1-q_{n-k}\Big[1-E\Big(\exp\Big\{-g\Big(\frac{Y_{n-k}}{n-k}\Big)\Big\}\Big|Y_{n-k}>0\Big)\Big]\Big)^{Z_k}\Big]\\
&=&E\Big[\exp\Big\{-\lim_{n\to\infty, k/n\to u}Z_kq_{n-k}\Big[1-E\Big(\exp\Big\{-g\Big(\frac{Y_{n-k}}{n-k}\Big)\Big\}\Big|Y_{n-k}>0\Big)\Big]\Big\}\Big].
\end{eqnarray*}
Then using \eqref{rate of ext prob with finit variance} and Lemma \ref{gamma limit}, we obtain
\begin{equation}\label{Laplace of point process}
\lim_{n\to\infty, k/n\to u}E[\exp\left\{-\langle g,V_{n,k}\rangle\right\}]=\lim_{n\to\infty, k/n\to u}E[L_{n,k}(g)]=E[\exp\{-\xi_u(1-L(g))\}],
\end{equation}
where $\xi_u$ is a random variable having Gamma distribution with parameters $\left(\dfrac{1-u}{ u},\gamma\right)$.
Then the Laplace transform of $\xi_u$ is given by (c.f. \cite[Example $2.15$]{Sato})
\[
L_{\xi_{u}}(\lambda)=Ee^{-\lambda\xi_u}=\Big(1+\dfrac{ u \lambda}{1-u}\Big)^{-\gamma},\qquad \lambda>0.
\]
Therefore,
\begin{eqnarray*}
&&E\left[\exp\left\{-\xi_u(1-L(g))\right\}\right]=\Big[1+\dfrac{u}{1-u}(1-L(g))\Big]^{-\gamma}
=(1-u)^{\gamma}[1-u L(g)]^{-\gamma}\\
&=&\sum_{k=0}^\infty \dfrac{(-\gamma)(-\gamma-1)\cdots (-\gamma-k+1)}{k!}(1-u)^{\gamma}(-u)^{k}L(g)^k\\
&=&Ee^{-\sum_{j=1}^{N^I_u}g(w_j)}=E\big[e^{-\langle g,V_u\rangle}\big].
\end{eqnarray*}
In conclusion, $V_{n,k}$ converges to $V_u$ in distribution as $n\to\infty,k/n\to u$. \qed

\bigskip

For $r>0$, put
 \begin{equation}\label{def-f-g}
 f(r)=r,\quad g_1(r)=r\wedge r^{-1}, \quad g_2(r)=1\wedge r^{2}.
 \end{equation}
\begin{remark}\label{from tild V to V}
Using the same argument as in the proof of Lemma \ref{Lemma of crm} for the random measure
$$
\widetilde{V}_{n,k}(\cdot):=\sum_{i=1}^{Z_k}{\rm I}_{\{Y_{n,k,i}>0\}}\Big(1\vee \Big(\dfrac{Y_{n,k,i}}{n-k}\Big)^2\Big) \delta_{\frac{Y_{n,k,i}}{n-k}}(\cdot),
$$
and using the fact that $h(r):=(1\vee r^2)$  is a continuous function on $[0,\infty)$,
we obtain that
$$
\widetilde{V}_{n,k}(dr)\stackrel{d}{\to}(1\vee r^2) V_u(dr)=: \widetilde{V}_u(dr)\quad \mbox{ in }\mathcal M.
$$
Since $g_1, g_2\in C_b^+[0,\infty)$,
$\langle g_1, \widetilde{V}_{n,k}\rangle=\langle f, V_{n,k}\rangle$
and $\langle g_2, \widetilde{V}_{n,k}\rangle=\langle f^2, V_{n,k}\rangle$,
we have
\begin{eqnarray}\label{conv for Vu}
&&(\langle f, V_{n,k}\rangle, \,\, \langle f^2, V_{n,k}\rangle )=(\langle g_1, \widetilde{V}_{n,k}\rangle, \,\, \langle g_2, \widetilde{V}_{n,k}\rangle )\nonumber\\
&\stackrel{d}{\to}& (\langle g_1, \widetilde{V}_{u}\rangle,\,\, \langle g_2, \widetilde{V}_{u}\rangle)=(\langle f, V_{u}\rangle,\,\, \langle f^2, V_{u}\rangle )
=\Big(\sum_{k=1}^{N^I_u}\omega_k ,\,\, \sum_{k=1}^{N^I_u}\omega_k^2\Big),
\end{eqnarray}
as $n\to\infty, k/n\to u$ with $u\in(0,1)$.
\end{remark}

Define the birth time $\tau_n$ of the oldest clan in
generation
$n$ by
\[
\tau_n=\inf\Big\{ 0\leq j\leq n; \sum_{l=1}^{I_j}Y_{n-j,l}^{(j)}>0\Big\}
\]
with the convention $\inf\emptyset=+\infty$.
The birth time of the oldest clan for stationary continuous state branching processes is studied in \cite[Corollary 4.2]{Chen Del}.
Using Lemma \ref{Lemma of crm}, it is easy to get the limit distribution of $\tau_n$.  Recall that $\gamma=2\beta/\sigma^2$.
\begin{corollary}\label{result for An}
Suppose Assumption \ref{offspring assumption} holds. We have
\[
\lim_{n\to\infty, k/n\to u}P(\tau_n>k)=P(N^I_u=0)=(1-u)^\gamma,\qquad 0<u<1.
\]
\end{corollary}
\noindent{\bf Proof.}
The event $\{\tau_n>k\}$ can be written as $\{V_{n,k}(1)=0\}$.  Thus
\[
\lim_{n\to\infty, k/n\to u}P(\tau_n>k)=\lim_{n\to\infty, k/n\to u}P(V_{n,k}(1)=0)=P(N^I_u=0)=(1-u)^{\gamma}.
\]
\qed

Define a function $w$ by
\begin{equation}\label{def-w}
w(r)=r\vee r^2,\quad r\in (0,\infty).
\end{equation}
We next consider the following random measures related to
immigrations after generation
 $k$,
$$
W_{n,k}(\cdot):=\sum_{j=k+1}^{n} \sum_{l=1}^{I_j}{\rm I}_{\big\{Y_{n-j, l}^{(j)}>0\big\}}w\Big(\dfrac{Y_{n-j, l}^{(j)}}{n-k}\Big)\delta_{\frac{Y_{n-j, l}^{(j)}}{n-k}}(\cdot), \qquad n>k.
$$
For each $(n,k)$ with $k<n$,
thanks to \eqref{decomp. at k},
we see that $W_{n,k}(\cdot)$ has the same distribution as the random measure
\begin{equation}\label{tild W}
\widetilde{W}_{n-k}(\cdot)
:=
\sum_{j=0}^{n-k-1} \sum_{l=1}^{I_j}{\rm I}_{\big\{Y_{j, l}>0\big\}}w\Big(\dfrac{Y_{j, l}}{n-k}\Big)
\delta_{\frac{Y_{j, l}}{n-k}}(\cdot),
\end{equation}
where $Y_{j, l},j\in\mathbb{N},l=1,2,\ldots$, are independent
and for each $j$, $Y_{j,l}, l=1,2,\ldots$, are identically distributed as $Y_j$, and where
$(Y_{j,l})_{j\geq 0, l\geq 1}$
are  independent of the immigration process
$(I_j)_{j\geq 0}$.
 By an argument very similar to that used in the proof of Lemma \ref{Lemma of crm}, we  get the following convergence in distribution result for the
random measures $(W_{n,k})_{n\geq k}$.

\begin{lemma}\label{lemma of immigration}
Suppose Assumption \ref{offspring assumption} holds.
Let $\zeta$ be the random measure defined  by
\[
\zeta(dr)= w(r)W(dr),
\]
where $W$ is a Poisson random measure with intensity $\frac{\gamma}{r}e^{-\frac{2r}{\sigma^2}}dr$
on $(0,\infty)$ and $w$ is the function defined in \eqref{def-w}.
Then $ W_{n,k}\stackrel{d}{\to} \zeta$ in $\mathcal M$  as  $n-k\to\infty$.
\end{lemma}
\noindent {\bf Proof.}
Since $W_{n,k}\stackrel{d}{=}\widetilde{W}_{n-k}$,
we have for
$g\in C_b^+[0,\infty)$,
\begin{eqnarray}\label{equal-Laplace}
E\big[\exp\big\{-\langle g,W_{n,k}\rangle\big\}\big]=E\Big[\exp\big\{-\langle g,\widetilde{W}_{n-k}\rangle\big\}\Big],
\end{eqnarray}
which means that we only need to consider the limit of the Laplace
functional
of $\widetilde{W}_{n}$ as $n\to\infty$.
For any  $g\in C^+_b[0,\infty)$,
put
$$
T_{n,j}(g)=E\Big[\exp\Big\{-w\left(\dfrac{Y_j}{n}\right)g\left(\dfrac{Y_j}{n}\right){\rm I}_{\{Y_j>0\}}\Big\}\Big],\quad j=0,1,\cdots, n-1.
$$
Then $0<T_{n,j}(g)<1$. By the definition \eqref{tild W} of $\widetilde{W}_{n}$,
\[
\exp\big\{-\langle g,\widetilde{W}_{n}\rangle\big\}=\exp\Big\{-\sum_{j=0}^{n-1} \sum_{l=1}^{I_j}w\left(\dfrac{Y_{j, l}}{n}\right)g
\left(\frac{Y_{j, l}}{n}\right){\rm I}_{\{Y_{j, l}>0\}}\Big\}.
\]
The Laplace transform of $\widetilde{W}_n$ can be written as
\begin{eqnarray}\label{Lap of immigration}
E\big[\exp\big\{-\langle g,\widetilde{W}_{n}\rangle\big\}\big]
=\prod_{j=0}^{n-1}E\big[T_{n,j}(g)^{I_j}\big]=\prod_{j=0}^{n-1}B\big(T_{n,j}(g)\big)=\exp\Big\{\sum_{j=0}^{n-1}\ln B\big(T_{n,j}(g)\big)\Big\},
\end{eqnarray}
where $B(s)=\sum_jb_j s^j, |s|<1$, is the probability generating function of $I_k, k\geq 0$.
We claim that
\begin{eqnarray}\label{conv 3}
\lim_{n\to\infty}\sum_{j=0}^{n-1}\ln B\big(T_{n,j}(g)\big)=\gamma\int_0^\infty\big(e^{-w(r)g(r)}-1\big)\dfrac{1}{r}e^{-\frac{2r}{\sigma^2}}dr.
\end{eqnarray}
Suppose for the moment the claim is true. Then by \eqref{Lap of immigration}, for any $ g\in C_b^+[0,\infty)$,
\begin{eqnarray*}
 \lim_{n\to\infty}E\big[\exp\big\{-\langle g,\widetilde{W}_{n}\rangle\big\}\big]=
\exp\Big\{\gamma\int_0^\infty\big(e^{-w(r)g(r)}-1\big)\frac{1}{r}e^{-\frac{2r}{\sigma^2}}dr\Big\}.
\end{eqnarray*}
And then using \eqref{equal-Laplace}, we have
\begin{eqnarray*}
\lim_{n-k\to\infty}E\big[\exp\big\{-\langle g,W_{n,k}\rangle\big\}\big]=
\exp\Big\{\gamma\int_0^\infty\big(e^{-w(r)g(r)}-1\big)\dfrac{1}{r}e^{-\frac{2r}{\sigma^2}}dr\Big\}.
\end{eqnarray*}
Since
$\int_0^\infty (w(r)\wedge 1) \frac{1}{r}e^{-\frac{2r}{\sigma^2}}dr<\infty$,
 it follows from \cite[Theorem 3.20]{Kal}
 that there is an infinitely divisible random measure $\zeta \in\mathcal M$
 represented as
$
\zeta(dr)=w(r)W(dr), r>0,
$
where $W$ is a Poisson random
measure with intensity  ${\rm I}_{\{r>0\}}\frac{\gamma}{r}e^{-\frac{2r}{\sigma^2}}dr$.  The Laplace functional of $\zeta$ is given by
\[
E\big[\exp\{-\langle g, \zeta\rangle\}\big]=\exp\Big\{\gamma\int_0^\infty \big(e^{-w(r)g(r)}-1\big)\dfrac{1}{r}e^{-\frac{2r}{\sigma^2}}dr\Big\},\quad \forall g\in C_b^+[0,\infty).
\]
Thus $W_{n,k}\stackrel{d}{\to}\zeta$
as $n-k\to\infty$.

Now we prove the claim \eqref{conv 3}.
By the mean value theorem, there exists $\xi_{n,j}\in (T_{n,j}(g),1)$ such that
\begin{eqnarray}\label{decom for B}
B\big(T_{n,j}(g)\big)-1&=&B'(\xi_{n,j})\big(T_{n,j}(g)-1\big)\nonumber\\
&=&\beta\big(T_{n,j}(g)-1\big)+(B'(\xi_{n,j})-\beta)\big(T_{n,j}(g)-1\big).
\end{eqnarray}
Thanks to the inequality $0<1-e^{-x}\leq x$ for $x>0$ and the fact that $\mbox{Var}(Y_j)=j\sigma^2$ (see \cite[Section $I.2$]{Ath Ney}), we have that for $0\leq j\leq n-1$,
\begin{eqnarray}\label{uniform bound}
0\leq 1-T_{n,j}(g)\leq \|g\|_\infty E\Big[w\big(\frac{Y_{j}}{n}\big)\Big]\leq \|g\|_\infty E\Big[\frac{Y_{j}}{n}+ \left(\frac{Y_{j}}{n}\right)^2\Big]\leq \frac{a\|g\|_\infty}{n},
\end{eqnarray}
for some  constant $a>0$. Thus $n(1-T_{n,j}(g))$ is bounded for $n>0$ and $j\leq n$.  Moreover
from \eqref{rate of ext prob with finit variance} and \eqref{limit of distri with finit var}, it follows that
for any $0<t<1$,
\begin{eqnarray*}
\lim_{n\to\infty} n[1-T_{n,[nt]}(g)]
&=& \lim_{n\to\infty} n P(Y_{[nt]}>0)E\Big[1-\exp\Big\{-w\big(\frac{Y_{[nt]}}{n}\big)g\big(\frac{Y_{[nt]}}{n}\big)\Big\}\Big|Y_{[nt]}>0 \Big]\\
&=&\dfrac{4}{(\sigma^2)^2 t}\int_0^\infty \left(1-e^{- w(rt)g(rt)}\right)e^{-\frac{2r}{\sigma^2 }}dr.
\end{eqnarray*}
Then by the dominated convergence theorem,
\begin{eqnarray}\label{conv 1}
\lim_{n\to\infty}\sum_{j=0}^{n-1}(T_{n,j}(g)-1)&=&\lim_{n\to\infty}\int_0^1 n(T_{n,[nt]}(g)-1)dt\nonumber\\
&=&\int_0^{1}\dfrac{4}{(\sigma^2)^2 t}dt\int_0^\infty \big(e^{- w(rt)g(rt)}-1\big)e^{-\frac{2r}{\sigma^2 }}dr\\
&=&\dfrac{2}{\sigma^2}\int_0^\infty\big(e^{-w(r)g(r)}-1\big)\dfrac{1}{r}e^{-\frac{2r}{\sigma^2}}dr.\nonumber
\end{eqnarray}
 Using \eqref{uniform bound} and the continuity of $B'(s)$
on $[0,1]$,
 we get that
 $B'(\xi_{n,j})-\beta$ converges to $0$ uniformly for $0\leq j\leq n$, as $n\to\infty$.  It has been shown in \eqref{conv 1} that $\sum_{j=0}^{n-1}\big|T_{n,j}(g)-1\big|$ converges. Therefore,
$\sum_{j=0}^{n-1}(B'(\xi_{n,j})-\beta)\big(T_{n,j}(g)-1\big)$ converges to $0$.  Thus, by \eqref{decom for B},
$\sum_{j=0}^{n-1}\big(B\big(T_{n,j}(g)\big)-1\big)$ and $\beta\sum_{j=0}^{n-1}\big(T_{n,j}(g)-1\big)$ have the same limit.  More precisely, from \eqref{conv 1}, it follows that
\begin{eqnarray*}
&&\lim_{n\to\infty}\sum_{j=0}^{n-1}(B\big(T_{n,j}(g)\big)-1)=\lim_{n\to\infty}\beta\sum_{j=0}^{n-1}\big(T_{n,j}(g)-1\big)\\
&=&\gamma\int_0^\infty\big(e^{-w(r)g(r)}-1\big)\dfrac{1}{r}e^{-\frac{2r}{\sigma^2}}dr.
\end{eqnarray*}
Meanwhile, since $-x\geq \ln (1-x)\geq -x-\dfrac{x^2}{1-x}$ for $0<x<1$, if
\begin{eqnarray}\label{error limit}
\lim_{n\to\infty}\sum_{j=0}^{n-1}\dfrac{\big[B\big(T_{n,j}(g)\big)-1\big]^2}{B\big(T_{n,j}(g)\big)}=0,
\end{eqnarray}
then
$\sum_{j=0}^{n-1}\ln B\big(T_{n,j}(g)\big)$ and $\beta\sum_{j=0}^{n-1}\big(T_{n,j}(g)-1\big)$ have the same limit  as $n\to\infty$,
and thus  the claim is true.
Now we prove \eqref{error limit}. By \eqref{uniform bound},
 for any $1/2<\delta<1$, there is $N>0$, such that for any $n>N, 0<j\leq n$,
$T_{n,j}(g)>\delta$. Since $B(s)$ is
an increasing continuous function on $[0,1]$
and $B(1)=1$, for any $\varepsilon>0$, we can choose $\delta$ above such that when $1>s>\delta$, $B(s)>1-\varepsilon$.  Therefore when $n>N$,
\[
0\leq \sum_{j=0}^{n-1}\dfrac{\big[B\big(T_{n,j}(g)\big)-1\big]^2}{B\big(T_{n,j}(g)\big)}\leq \dfrac{\varepsilon}{B(\frac{1}{2})}\sum_{j=0}^{n-1}\big[1-B\big(T_{n,j}(g)\big)\big].
\]
Then \eqref{error limit} follows from the convergence of $\sum_{j=0}^{n-1}\big[1-B\left(T_{n,j}(g)\right)\big]$ and the arbitrariness of $\varepsilon$.
\qed

\begin{remark}\label{conv for Wu}
$(1)$  Let $\tilde g_1(r)=1\wedge r^{-1}, r>0$ and $\tilde g_2(r)=1\wedge r, r>0$. Then
$\tilde g_1, \tilde g_2\in C_b^+[0,\infty)$.
Thanks to Lemma \ref{lemma of immigration} and the facts
$\tilde g_1(r)w(r)=r=f(r)$ and $\tilde g_2(r)w(r)=r^2=f^2(r)$ for $r>0$, we get that
$$
\big(\langle \tilde g_1,W_{n,k}\rangle, \, \langle \tilde g_2,W_{n,k}\rangle\big)\stackrel{d}{\to}\big(\langle \tilde g_1,\zeta\rangle, \, \langle \tilde g_2, \zeta\rangle\big)
=\big(\langle f,W\rangle, \, \langle f^2, W\rangle\big),\quad {\rm as}\,\, n-k\to\infty.
$$
$(2)$  We observe that
$$
\dfrac{n-k}{n}\big[\langle f, V_{n,k}\rangle+
\langle \tilde g_1, W_{n,k}\rangle
\big]=\dfrac{Z_n}{n},
$$
where $f$ and $\tilde g_1$ are defined as above.
Since $ V_{n,k}$ and $W_{n,k}$ are independent, from Lemma \ref{Lemma of crm} and Lemma \ref{lemma of immigration}, it follows that for any $\lambda>0$.
\begin{eqnarray*}
&&\lim_{n\to\infty}E\Big[\exp\Big\{-\lambda\frac{Z_n}{n}\Big\}\Big]=\lim_{n\to\infty, k/n\to u}E\exp\Big\{-\lambda\dfrac{n-k}{n}(\langle f, V_{n,k}\rangle+\langle \tilde g_1, W_{n,k}\rangle)\Big\}\\
&=&\lim_{n\to\infty, k/n\to u}E\Big[\exp\big\{-\lambda\dfrac{n-k}{n}\langle f, V_{n,k}\rangle\big\}\Big]\lim_{n\to\infty, k/n\to u}E\Big[\exp\big\{-\lambda\dfrac{n-k}{n}\langle \tilde g_1, W_{n,k}\rangle\big\}\Big]\\
&=&E[\exp\{-\lambda(1-u)\langle f, V_u\rangle\}]E[\exp\{-\lambda(1-u)\langle f, W\rangle\}]\\
&=&\Big(\dfrac{\lambda+\frac{2}{\sigma^2}}{\lambda(1-u)+\frac{2}{\sigma^2}}\Big)^{-\gamma}
\Big(\dfrac{\lambda(1-u)+\frac{2}{\sigma^2}}{\frac{2}{\sigma^2}}\Big)^{-\gamma}
=\Big(1+\dfrac{\lambda}{\frac{2}{\sigma^2}}\Big)^{-\gamma},
\end{eqnarray*}
where the last term
is the Laplace transform of the Gamma distribution with parameters $(\frac{2}{\sigma^2},\gamma)$. This is consistent with Lemma \ref{gamma limit}.
\end{remark}

\section{Proofs of the main results }

{\bf Proof of Theorem \ref{main limit thm1}: }
Let $f$ be the function defined in \eqref{def-f-g}, and let $\tilde{g}_1, \tilde{g}_2$ be the functions defined in Remark \ref{conv for Wu}$(1)$.
The random variable in Lemma \ref{main thm 1} can be expressed in terms of
the  random measures defined in Section $2$, and then we have
\begin{eqnarray*}
\dfrac{\sum_{i=1}^{Z_k}\big(Y_{n,k,i}\big)_2+\sum_{j=1+k}^{n-1}\sum_{l=1}^{I_j}\big(Y_{n-j,l}^{(j)}\big)_2}{(Z_n)_2}
=
\dfrac{\langle f^2, V_{n,k}\rangle-\frac{1}{n-k}\langle f, V_{n,k}\rangle+\langle \tilde{g}_2, W_{n,k}\rangle-\frac{1}{n-k}\langle \tilde{g}_1, W_{n,k}\rangle}
{\big[\langle f, V_{n,k}\rangle+\langle \tilde{g}_1, W_{n,k}\rangle\big]^2-\frac{1}{n-k}\big[\langle f, V_{n,k}\rangle+\langle \tilde{g}_1, W_{n,k}\rangle\big]}.
\end{eqnarray*}
Since $(V_{n,k})_{n>k}$ and $(W_{n,k})_{n>k}$ are independent and  $0<\dfrac{\sum_{i=1}^{Z_k}\big(Y_{n,k,i}\big)_2+\sum_{j=1+k}^n\sum_{l=1}^{I_j}\big(Y_{n-j,l}^{(I_j)}\big)_2}{(Z_n)_2}\leq 1$ is a bounded continuous function of $(\langle f,V_{n,k}\rangle, \langle f^2,V_{n,k}\rangle, \langle \tilde{g}_1, W_{n,k}\rangle, \langle \tilde{g}_2, W_{n,k}\rangle)$, according to
Remark \ref{from tild V to V} and Remark \ref{conv for Wu}, for $u\in (0,1)$,
\begin{eqnarray*}
&&\lim_{n\to\infty, k/n\to u}\dfrac{\langle f^2, V_{n,k}\rangle-\frac{1}{n-k}\langle f, V_{n,k}\rangle+\langle \tilde{g}_2, W_{n,k}\rangle-\frac{1}{n-k}\langle \tilde{g}_1, W_{n,k}\rangle}{\big[\langle f, V_{n,k}\rangle+\langle \tilde{g}_1, W_{n,k}\rangle\big]^2-\frac{1}{n-k}\big[\langle f, V_{n,k}\rangle+\langle \tilde{g}_1, W_{n,k}\rangle\big]}=\dfrac{\langle f^2, V_u\rangle+\langle f^2, W\rangle}{\big[\langle f, V_{u}\rangle+\langle f, W\rangle\big]^2}
\end{eqnarray*}
in distribution.
It follows from  Lemma \ref{gamma limit} that $\lim_{n\to\infty}P(Z_n>1)=1$.
The results of this theorem follow from Lemma \ref{main thm 1}.
\bigskip

\noindent {\bf Proof of Theorem \ref{result for taun}:}
If all the particles in generation $n$ have  the same ancestor,
then they must be
descendants of one immigrant before generation
$n$.
Thus
\begin{eqnarray*}
\big\{A^I_n<\infty, Z_n>0\big\}
 \subset \big\{Y_{n-j,l}^{(j)}=0\, \mbox{for all but one pair } (j,l), 0\leq j\leq n, 1\leq l\leq I_j\big\}.
\end{eqnarray*}
Then we only need to prove that the probability of the event on the right hand side
converges to $0$.
Recall that $q_n=P(Y_n>0)$.  Set
$a_n=1-q_n=P(Y_n=0)$.
 Then
\begin{eqnarray}\label{upp}
&&P\Big(Y_{n-j,l}^{(j)}=0\, \mbox{for all but one pair} (j,l), 0\leq j\leq n, 1\leq l\leq I_j\Big)\nonumber\\
&=& E\Big[\sum_{j=0}^{n}\prod_{k\neq j}P(Y_{n-k}=0)^{I_k}I_jP(Y_{n-j}=0)^{I_j-1}P(Y_{n-j}>0)\Big]\nonumber\\
&=&\Big[\prod_{k=0}^{n}B(a_k)\Big]\Big[\sum_{j=0}^{n}\dfrac{B'(a_j)}{B(a_j)}q_j\Big],
\end{eqnarray}
where $B(a_0)=B(0)=b_0$ and $B'(a_0)=B'(0)=b_1$.
From \eqref{rate of ext prob with finit variance}, we know $q_k=1-a_k\sim \frac{2}{\sigma^2 k}$ as $k\to\infty$.  In addition, since $B(s)=1+\beta(s-1)+
o(1-s)$ as $s\to 1-$,
\begin{eqnarray}\label{limit of B}
\lim_{j\to\infty}j(1-B(a_j))=\lim_{j\to\infty}\beta j(1-a_j)+ o(j(1-a_j))=\gamma>0.
\end{eqnarray}
Therefore, there exists some $N\in\mathbb N$, such that when $k\geq N$, $k(1-B(a_k))>\gamma/2$, which implies that
$B(a_k)<1-\frac{\gamma}{2k}$ for $k\geq N$.  Noticing that $B(a_k)\leq 1$,
the first factor on the right-hand side of \eqref{upp}  can be estimated as follows:
\[
\prod_{j=0}^{n}B(a_j)\leq \prod_{j=N}^{n}B(a_j)\leq \prod_{j=N}^{n}\big(1-\dfrac{\gamma}{2j}\big)
=\exp\Big\{\sum_{j=N}^{n}\ln(1-\dfrac{\gamma}{2j})\Big\},\quad n>N.
\]
Since $\ln(1-x)<-x$ for $0<x<1$, we have
\[
\sum_{k=N}^{n}\ln(1-\dfrac{\gamma}{2k})\leq -\sum_{k=N}^{n}\dfrac{\gamma}{2k}\leq -L(\ln n- \ln N),
\]
for some constant $L>0$.  As a result, there exists $C_1>0$, such that
\begin{eqnarray}\label{term 1}
\prod_{k=0}^{n}B(a_k)\leq C_1 \cdot n^{-L}.
\end{eqnarray}
Since $a_k$ is nondecreasing in $k$
and converges to $1$ as $k\to \infty$, and $B'(s)$ is a continuous function on $[0,1]$,
\[
\lim_{j\to\infty}\dfrac{B'(a_j)}{B(a_j)}=B'(1)=\beta.
\]
The the second factor on the right-hand side of \eqref{upp} has the following upper bound:
\begin{eqnarray}\label{term 2}
\sum_{j=1}^{n}\dfrac{B'(a_j)}{B(a_j)}q_j\leq C_2\sum_{j=1}^nq_j\leq C_3 \sum_{j=1}^n\dfrac{1}{j}\leq C_3(1+\ln n),
\end{eqnarray}
for some positive constants $C_2$ and $C_3$.
Combining \eqref{term 1} and \eqref{term 2}, we obtain
\[
\lim_{n\to\infty}P\Big(Y_{n-j,l}^{(j)}=0\, \mbox{for all but one pair } (j,l),\, 0\leq j\leq n, 1\leq l\leq I_j\Big)=0.
\]
We finish the proof.\qed

\vspace{.1in}
	\textbf{Acknowledgment}:
We thank the referees for very helpful comments and suggestions.

\end{document}